\documentclass[12pt,a4paper,oneside]{amsart}

\ifdefined\SMART
\usepackage[paperwidth=12cm,paperheight=24cm,left=5mm,right=5mm,top=10mm,bottom=5mm]{geometry}
\else
\calclayout
\fi

\usepackage{amssymb,latexsym,amsfonts,textcomp,fancyhdr,calc,graphicx}
\usepackage{pdfpages,amsfonts,amsthm,amssymb,latexsym,graphicx,leftidx,fancyhdr}
\usepackage{amsmath,amstext,amsthm,amssymb,amsxtra,mathabx,xcolor}
\usepackage[allcolors=blue,allbordercolors=blue,pdfborderstyle={/S/U/W 1}]{hyperref}
\usepackage[ocgcolorlinks]{ocgx2}
\usepackage{dsfont}
\usepackage[framemethod=tikz]{mdframed}
\usepackage{asymptote}
\usepackage{enumerate}

\usepackage[normalem]{ulem}
\usepackage{cite}
\usepackage{comment}

\usepackage{txfonts,pxfonts,tikz} 
\usepackage{tgschola}
\usepackage[T1]{fontenc}

\usepackage{mathtools}

\usepackage{bbm}

\newtheorem{theorem}{Theorem}
\newtheorem{corollary}{Corollary}
\newtheorem{lemma}{Lemma}
\newtheorem{example}{Example}
\newtheorem{question}{Question}
\newtheorem{conjecture}{Conjecture}
\newtheorem{proposition}{Proposition}
\theoremstyle{definition}

\newcommand{\beql}[1]{\begin{equation}\label{#1}}
\newcommand{\eeq}{\end{equation}}
\renewcommand{\comment}[1]{}

\newcommand{\Abs}[1]{{\left|{#1}\right|}}

\newcommand{\Linf}[1]{{\left\|{#1}\right\|_\infty}}

\newcommand{\Set}[1]{{\left\{{#1}\right\}}}

\newcommand{\RR}{{\mathbb R}}

\newcommand{\ZZ}{{\mathbb Z}}

\newcommand{\one}{{\mathbbm{1}}}

\newcommand{\dens}{{\rm dens\,}}
\newcommand{\supp}{{\rm supp\,}}

\newcommand{\ft}[1]{\widehat{#1}}

\newcommand{\tf}[1]{\widecheck{#1}}

\newcommand{\wt}[1]{\widetilde{#1}}

\renewcommand{\Re}{{\rm Re\,}}

\newcounter{rem}
\newcounter{step}
\newcounter{othm}
\def\theothm{\Alph{othm}}

\newcounter{mysec}
\newcounter{mysubsec}[mysec]
\begin{document}

\sloppy

\title[{Non-spectrality of some curves and lines segments }]{Non-spectrality of some piecewise smooth curves and unions of line segments }

\author{Mihail N. Kolountzakis}
\address{\href{http://math.uoc.gr/en/index.html}{Department of Mathematics and Applied Mathematics}, University of Crete,\\Voutes Campus, 70013 Heraklion, Greece,\newline and \newline \href{https://ics.forth.gr/}{Institute of Computer Science}, Foundation of Research and Technology Hellas, N. Plastira 100, Vassilika Vouton, 700 13, Heraklion, Greece}
\email{kolount@uoc.gr}

\author{Chun-Kit Lai}
\address{\href{https://math.sfsu.edu/}{Department of Mathematics}, San Francisco State University, San Francisco,CA 94132.}
\email{cklai@sfsu.edu}

\makeatletter
\@namedef{subjclassname@2020}{\textup{2020} Mathematics Subject Classification}
\makeatother
\subjclass[2020]{42C15, 42C30}

\keywords{Spectrality, piecewise smooth curves, tempered distributions}

\begin{abstract}
We develop a systematic study about the spectrality of measures supported on piecewise smooth curves by studying the support of the tempered distributions arising from the tiling equation of some singular spectral measures. In doing so, we show that the arc-length measures of all closed polygonal lines are not spectral. {In particular, the boundary of a square is not spectral. We also show that the ``plus space'' (two crossing line segments) is not spectral.} Furthermore, our theory also shows that the arc length measures on {smooth} convex curves with finitely many transverse self-intersections are not spectral. Finally, several natural open questions about the spectrality of singular measures and {piecewise} smooth curves will also be discussed. 
\end{abstract}

\date{\today}

\maketitle


\section{Introduction}\label{s:intro}

In this paper, a measure $\mu$ is always referring to a compactly supported probability measure in $\RR^d$.  We say that $\mu$ is a {\bf spectral measure} if there exists a countable set $\Lambda\subset \RR^d$ such that $E(\Lambda): = \{e^{2\pi i \lambda \cdot x}: \lambda\in\Lambda\}$ forms an orthonormal basis for $L^2(\mu)$. A measurable set $\Omega$ with finite positive Lebesgue measure is a {\bf spectral set} if the normalized Lebesgue measure on $\Omega$ is a spectral measure. The study of spectral sets of was first initiated by Fuglede in 1974 \cite{fuglede1974operators}, who famously proposed the conjecture that {\it spectral sets are equivalent to translational tiles on $\RR^d$}. The conjecture was disproved in both directions on $\RR^d$ with $d\ge 3$ \cite{tao2004fuglede,kolountzakis2006tiles,kolountzakis2006hadamard,farkas2006onfuglede,farkas2006tiles}, but the research effort to explore the exact geometric relationships between spectral sets and translational tiling has not declined. {In one recent major result,} Lev and Matolcsi \cite{lev2022fuglede} showed that Fuglede's conjecture holds for convex domains.  Readers can refer to the survey \cite{kolountzakis2024orthogonal} for different aspects about the Fuglede's conjecture.  

The generalized concept of spectral measures gained its popularity when Jorgensen and Pedersen discovered that {the  self-similar Cantor measure generated by dividing $[0,1]$ into 4 sub-intervals of equal length and keeps the first and third one is a spectral measure}, while the standard middle-third Cantor measure is not a spectral measure in 1998 \cite{jorgensen1998dense}. This work was followed up by Strichartz \cite{strichartz2000mock} and {\L}aba and Wang \cite{laba2002spectral-cantor}.  There {is now a vast} literature in this direction.  Readers can refer to the survey \cite{DLW17} for more recent advances about fractal spectral measures.  The connection with tiling appears to be weaker, but it exists. Gabardo and Lai \cite{GL14} showed that if the pair of measures  $(\mu,\nu)$ ``tiles" $[0,1]{^d}$ in the sense that $\mu\ast \nu = {\one_{[0,1]^d}\,dx}$, then both $\mu$ and $\nu$ are spectral measures. In the same paper, they also proposed a version of the generalized Fuglede's conjecture for singular spectral measures. A verification of {this conjecture in} a special case of Cantor-Moran measures was recently given in \cite{ALZ25}.

The purpose of this paper is to study the spectrality of singular measures supported on piecewise smooth curves.  This first appeared in Lev's paper \cite{lev2018fourier}.  Iosevich, Lai, Liu and Wyman showed that boundary of the unit disk {never admits a Fourier frame  (generalization of the orthonormal basis of exponentials), hence it is non-spectral}\cite{iosevich2022fourier}. {In the same paper}, they also showed that the boundary of any polygons always admit a Fourier frame. This naturally {led to the} question to determine if the boundary of a square, equipped with the 1-dimensional Lebesgue measure, is a spectral measure. {Notice a well-known fact that the square domain is spectral, thus it becomes an interesting question to determine if the spectrality of the domain can induce a spectrum on its boundary. We will answer these questions in this article in the negative.} 

\begin{figure}[h]
\ifdefined\SMART\resizebox{4cm}{!}{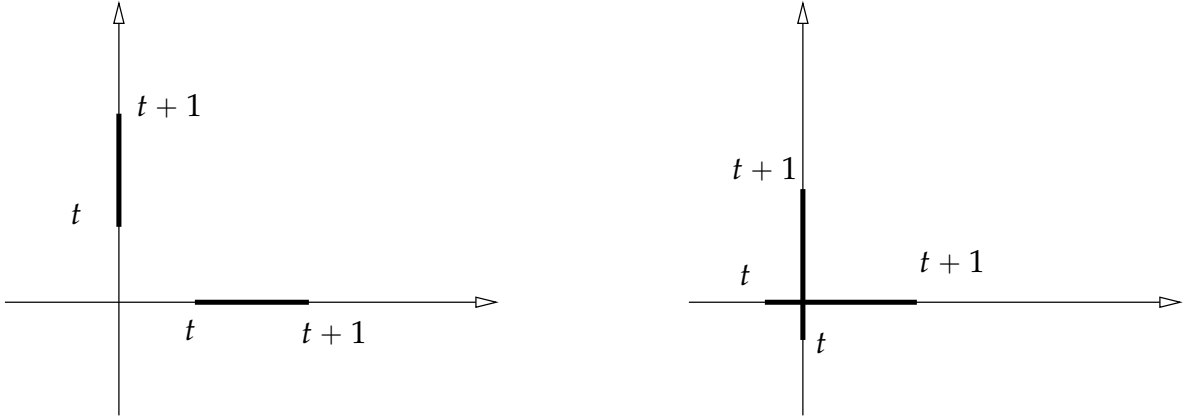}\else
\input sticks.pdf_t
\fi
\caption{The measure that is arc-length on two equal-length line-segments. Some of the measures on the left are spectral, depending one where the starting points of the segments are compared to their length. In the crossing case on the right they are never spectral. The symmetric crossing case, $t=-1/2$, is called the ``plus-space''.}\label{fig:sticks}
\end{figure}

Lai, Liu and Prince \cite{lai2021spectral} initiated a new study about the spectrality of symmetric additive spaces with measure supported on two orthogonal lines defined 
\begin{equation}\label{eq-symmetric-measure}
\mu = \frac12 (\one_{[t,t+1]}~dx\times \delta_0+ \delta_0\times \one_{[t,t+1]}~dx).    
\end{equation}
They proved that $\mu$ is a spectral measure when $t = 0$. However, the case where the two lines overlap, i.e. $-1/2\le t<0$, is much more subtle. With a series of work{s} by Ai-Lu-Zhou \cite{ai2023spectrality}, Kolountzakis and Wu \cite{kolountzakis2025spectralitymeasureconsistingline} and recently by Lu \cite{lu2025frameboundspectralgap}, it is now shown that all are not spectral when $-1/2\le t<0$. See Fig.\ \ref{fig:sticks}.

Apart from {spectra}, our results also deal with tight frames of exponentials, a slightly more general notion of basis.  We say that $\mu$ is a {\bf tight-frame spectral measure} if there exists a countable set $\Lambda\subset \RR^d$, which we call a {\bf tight-frame spectrum}, such that $E(\Lambda): = \{e^{2\pi i \lambda \cdot x}: \lambda\in\Lambda\}$ forms a tight frame for $L^2(\mu)$ in the sense that there exists a constant $A>0$ such that 
\begin{equation}\label{tight-frame}
\sum_{\lambda\in \Lambda} \left|\int f(x)e^{-2\pi i \lambda\cdot x}~d\mu(x)\right|^2  = A \int|f(x)|^2~d\mu(x), \ \forall f\in L^2(\mu). 
\end{equation}
Clearly an orthonormal basis is a tight frame. But the converse is not true. For example, we can take union of two orthonormal basis to form a tight frame of constant $2.$ For general frame theory, readers can refer to \cite{Heilbook}. For tight frame spectral measures, general theory was established in \cite{HLL2013} and \cite{dutkay2014uniformity}. In \cite{dutkay2014uniformity}, it was shown that if $\mu$ is tight-frame spectral and absolutely continuous with respect to the Lebesgue measure, then $\mu$ must be a constant multiple of the Lebesgue measure on its support. The consideration of tight-frame spectrality allows us to deal with spectrality of curves by focusing on a subset of {the curve} (see Corollary \ref{corollary1} below).

We say that a finite union of line segments is {(tight-frame)} spectral if the naturally induced one-dimensional Lebesgue measure is a  
{(tight frame)} spectral measure. Our main result in this paper is to settle the spectrality questions we mentioned. {Moreover, we determine non-spectralities for much more general classes of line segment collections.}

\begin{theorem}\label{theorem1}
\begin{enumerate}
    \item A finite union of line segments that forms a closed curve{, self-intersecting or not,} cannot be  tight-frame spectral.
    \item A finite union of line segments containing three lines that starts at the same point and points in distinct directions cannot be tight-frame spectral.
\end{enumerate}
\end{theorem}

\begin{figure}[h]
\ifdefined\SMART\resizebox{4cm}{!}{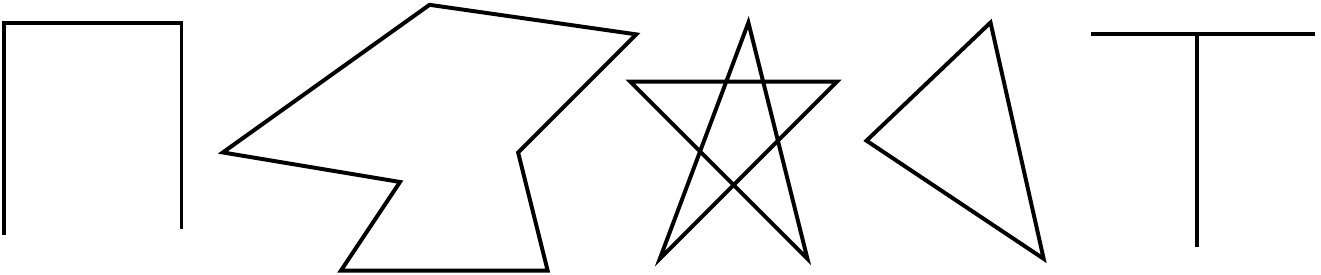}\else
\resizebox{15cm}{!}{\input polygons.pdf_t}
\fi
\caption{Some examples of non-spectral collections of line segments. All but the first one from the left, the ``$\Pi$-shape'', are covered by Theorem \ref{theorem1}. The $\Pi$-shape is explained in Section \ref{ss:two-cases}}\label{fig:polygons}
\end{figure}

In the second part, we remark that the theorem also {includes cases in which} two of the three vectors are in opposite directions to each other. Hence, it shows that the ``plus space'' corresponding to $t = -1/2$  {(refer to Fig.\ \ref{fig:sticks})} cannot be tight-frame spectral, providing another independent proof after Lu's proof \cite{lu2025frameboundspectralgap}.  
Some more examples following from our Theorem \ref{theorem1} are shown in Fig.\ \ref{fig:polygons}.

Clearly, the first part of the theorem showed that the boundary of a square cannot admit a tight frame of exponentials. Closed polygonal curves can be generalized to certain smooth curve{s} of positive curvature.

\begin{theorem}\label{theorem2}
Let $\gamma$ be a smooth closed planar curve with positive curvature such that there are only finitely many self-intersections {all of which are transverse}. Then the induced arc-length measure on $\gamma$ is not a {tight-frame} spectral measure. 
\end{theorem}

\begin{figure}[h]
\ifdefined\SMART\resizebox{4cm}{!}{\input closed-curve.pdf_t}\else
\resizebox{4cm}{!}{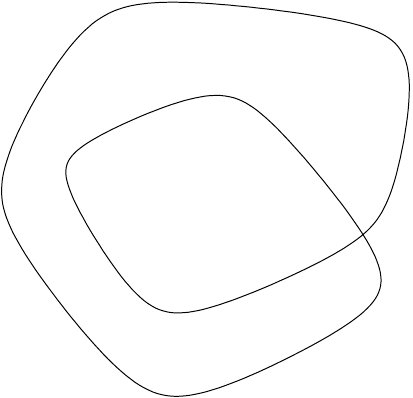}
\fi
\caption{A smooth closed curve in the plane with positive curvature everywhere and a transverse self-intersection.}\label{fig:closed-curve}
\end{figure}

Transverse intersections mean that the tangent vectors from different time{s visiting} the intersection {point} are not parallel to each others ({s}ee Section \ref{Section-curves} for the precise definition).  This result is slightly more general than the results in \cite{iosevich2022fourier} in the sense that {the curves we deal with} can be self-intersecting, such as the one in Fig.\ \ref{fig:closed-curve}.

Finally, we observe a simple property about tight frame spectral measures. Suppose that $E(\Lambda) = \{e^{2\pi i \lambda\cdot x}: \lambda\in\Lambda\}$ is a tight frame for the measure $\frac12(\mu+\nu)$ {and assume that $\supp\mu \cap \supp\nu$ has zero $\mu$-measure and $\nu$-measure}. By restricting only to functions on $L^2(\mu)$ {(i.e. extended by 0 off $\supp\mu$)}, we immediately see that $E(\Lambda)$ is also a tight frame spectrum for $L^2(\mu)$ with another positive constant.  Because of this, the following corollary is immediate:

\begin{corollary}\label{corollary1}
Let $\mu$ be the arc length measure supported on a finite union of smooth curves (not necessarily connected). Suppose that the union contains curves described in Theorem \ref{theorem1} and \ref{theorem2}. Then $\mu$ is not tight-frame spectral.
\end{corollary}

\subsection{Idea of the proof and function tilings.} The proof of our theorems requires a generalization of a well-known fact in the tiling theory of integrable functions to a non-integrable setting. 
Let us denote $\wt\mu$ the measure $\wt\mu (E) = {\mu (-E)}$  and define the Fourier transform of the measure $\mu$  to be 
$$
\widehat{\mu}(\xi)  = \int e^{-2\pi i \xi \cdot x}d\mu(x). 
$$
The convolution between two measures $\mu$ and $\nu$ is given by 
$$
\mu\ast\nu (E) = \int \one_E (x+y) d\mu (x)d\nu (y).
$$  
Recall that if $f\in L^1(\RR^d)$ and $\Lambda$ is a countable discrete set such that 
\begin{equation}\label{general-tiling}
\delta_{\Lambda}\ast f(x) = \sum_{\lambda\in \Lambda} f(x-\lambda) = w ~~\text{a.e.}
\end{equation}
for some $w$, then we can conclude that $\supp\ft{\delta_{\Lambda}}$ with $\ft{\delta_{\Lambda}}$ regarded as a tempered distribution, satisfies
$$
\supp\ft{\delta_{\Lambda}} \subset \{0\}\cup \{\xi\in\RR^d: \widehat{f}(\xi) = 0\}. 
$$
This can be formally deduced by taking {the Fourier transform of both sides of \eqref{general-tiling} to obtain} $\ft{\delta_{\Lambda}} \ft{f} = w \delta_0$. A rigorous proof of this fact can be found in \cite[Theorem 4.1]{kolountzakis2016non}.
{Recall that \cite{jorgensen1998dense}} $\mu$ is a spectral measure with a spectrum $\Lambda$ if and only if
\begin{equation}\label{conv-eq}
\delta_{\Lambda}\ast|\ft{\mu}|^2 = 1.
\end{equation}
{If $\mu$ admits a tight frame $\{e^{2\pi i \lambda\cdot x}: \lambda\in \Lambda\}$, by taking $f(x) = e^{2\pi i \xi\cdot x}$ into (\ref{tight-frame}), we obtain also a tiling equation
\begin{equation}\label{conv-eq1}
\delta_{\Lambda}\ast|\ft{\mu}|^2 = A.
\end{equation}}
Notice that the distributional Fourier transform of $|\ft{\mu}|^2$ is equal to the measure $\mu\ast\wt{\mu}$. By formally taking Fourier transform to the convolutional equation {\eqref{conv-eq}} and (\ref{conv-eq1}), the following conjecture may be true:

\begin{conjecture}\label{conjecture1}
Let $\Lambda$ be a {tight-frame} spectrum for a measure $\mu$ in $\RR^d$. Then 
$$
\supp\ft{\delta_{\Lambda}} \subset \{0\}\cup (\supp (\mu\ast\wt{\mu}))^C. 
$$
\end{conjecture}
{Above}, the support of a measure $\mu$ is the unique closed support of $\mu$, which is the set of points $x$ such that $\mu (B(x,r))>0$ for all $r>0$.  When $\mu$ is an absolutely continuous measure {with respect to Lebesgue measure, \cite{dutkay2014uniformity} showed that the density must be constant function, so it must be in  $L^2(\RR^d)$. It follows that} $|\ft{\mu}|^2$ is integrable and this conjecture is correct. We are unable to prove this conjecture in general. Instead, we have the following weaker version of the conjecture: 

\begin{theorem}\label{prop-non-spectral}
 Suppose $\mu$ is a {tight-frame} spectral measure on $\RR^d$ with a {tight-frame} spectrum $\Lambda$ {such that} $\mu*\wt\mu$ is absolutely continuous in the open set $U$ and has a smooth, strictly positive density therein.  Then $
\supp\ft\delta_\Lambda \cap U \subseteq \Set{0}.$    
\end{theorem}\label{lm:support}

{By a} {\bf spectral gap} $a>0$ for a tempered distribution $T$, we mean that $\ft{T}$ vanishes on a {punctured} open ball $B(0,a)\setminus \{0\}$. If a spectral gap does not exist, then we say that $T$ has a zero spectral gap.  

\begin{theorem}\label{prop-spectral-gap}
Let $\Lambda$ be a countably discrete set and $\Lambda$ {be} a {tight-frame} spectrum for {a} singular measure $\mu$. Then the spectral gap of $\delta_\Lambda $ is zero. 
\end{theorem}

 Conjecture \ref{conjecture1} and Theorem \ref{prop-spectral-gap} suggest a strategy of showing that a singular measure is non-spectral or non-tight-frame spectral. 
 
\begin{proposition}\label{prop-spectral-gap-1} Let $\mu$ be a singular measure. Suppose that Conjecture \ref{conjecture1} holds and  the support of $\mu\ast\wt\mu$ covers a neighborhood of the origin. Then $\mu$ cannot be tight-frame spectral.
\end{proposition}

 \begin{proof}
Suppose that $\Lambda$ is a tight-frame spectrum for $\mu$. From Conjecture \ref{conjecture1} and the assumption of the support, a spectral gap for $\delta_{\Lambda}$ exists. However, this is a contradiction to Theorem \ref{prop-spectral-gap} since $\mu$ is a singular measure.       
 \end{proof}

Despite Conjecture \ref{conjecture1} being open, we are able to adopt the same strategy to prove Theorem \ref{theorem2} via Theorem \ref{prop-non-spectral}. We will show that under the assumption of Theorem \ref{theorem2}, $\mu\ast\wt\mu$ has a smooth density around a punctured neighborhood of the origin (see Theorem \ref{theorem-smooth}).  To prove Theorem \ref{theorem1}, we will also adopt the same strategy. However, we will see that the measure $\mu\ast\wt{\mu}$ contains a singular part supported on some lines, and we will need a characterization of distributions supported on the lines to handle this situation (see Subsection \ref{subsection-line-distribution}).

We organize the paper as follows: In Section \ref{s:tool}, we will lay out our tools and previous results required for our study. In Section \ref{s:thm34}, we will prove Theorem 3 and Theorem 4. These two sections laid down the foundation to prove Theorem \ref{theorem1} in Section \ref{s:line} and Theorem \ref{theorem2} in Section \ref{Section-curves}. Finally, we will discuss some open problems in Section \ref{s:problems}

\section{Preliminaries}\label{s:tool}
We will provide all the basic tools required in this paper. In particular, we will provide the background on tempered distributions required for the study of the spectral measures, the properties about distributions supported on 1-dimensional subspaces and finally a change of variable formula. 

\subsection{Densities of countably discrete sets.}   For a countably discrete set $\Lambda\subset \RR^d$, the upper and lower {\bf Beurling densities} {are} defined to be 
$$
D^+(\Lambda) = \limsup_{r\to\infty} \sup_{x\in\RR^d} \frac{\#(\Lambda\cap B(x,r))}{r^d}, \ D^-(\Lambda) = \liminf_{r\to\infty} \inf_{x\in\RR^d} \frac{\#(\Lambda\cap B(x,r))}{r^d}.
$$
 $\Lambda$ has a {\bf uniformly bounded density} if dens$(\Lambda) =D^+(\Lambda) = D^-(\Lambda)<\infty$. $\Lambda$ is called {\bf translationally-bounded} if
 $$\sup_{x\in\RR^d}\#(\Lambda\cap B(x,1))<\infty.$$
 It is known that $\Lambda$ is translationally-bounded if and only if $D^+(\Lambda)<\infty$. 
 We also recall that $\Lambda$ is called {\bf separated} if there exists $\delta>0$ such that $|\lambda-\lambda'|\ge \delta$ for all distinct $\lambda,\lambda'\in\Lambda$. We begin with a simple lemma.

 \begin{lemma}\label{lemma-schwartz-bounded}
     Let $\psi$ be a Schwartz function and let $\Lambda$ be translationally bounded. Then 
     $$\delta_{\Lambda}\ast \psi (x) = \sum_{\lambda\in\Lambda}\psi (x-\lambda)
     $$
     is a bounded function in $\RR^d$. 
 \end{lemma}

 \begin{proof}Note that $D^{+}(\Lambda)<\infty$. {Then}
 $$
 \#(\Lambda\cap B(x,2^n)) \lesssim  2^{nd}
 $$
where the implicit constant is independent of $n$ {and $x$}.  As $\psi$ has a rapid decay, $|\psi(x)|\lesssim  (1+|x|)^{-(100d)}$.  Hence, 
$$    \sum_{\lambda\in\Lambda}|\psi (x-\lambda)| = \sum_{\lambda\in B(x,1)}|\psi(x-\lambda)|+\sum_{n=1}^{\infty}\sum_{\lambda\in B(x,2^n)\setminus B(x,2^{n-1})}|\psi(x-\lambda)|
$$
The first term is clearly bounded independent of $x$ since $\sup_{x\in\RR^d}\#(\Lambda\cap B(x,1))<\infty$. For the second sum, 
$$
\sum_{n=1}^{\infty}\sum_{\lambda\in B(x,2^n)\setminus B(x,2^{n-1})}|\psi(x-\lambda)|\lesssim  \sum_{n=1}^{\infty} \frac{2^{nd}}{(1+2^{n-1})^{100d}}\lesssim  \sum_{n=1}^{\infty} \frac{1}{2^{99dn}}<\infty.
$$
This completes the proof. 
\end{proof}

For $\Lambda \subseteq \RR^d$ of bounded density we write
$$
\delta_\Lambda = \sum_{\lambda\in\Lambda} \delta_\lambda.
$$
This is a tempered distribution which is locally a measure and we {can} therefore speak of its Fourier Transform, $\ft{\delta_\Lambda}$, which is also a tempered distribution, but we notice that it is not necessarily locally a measure. Readers can refer to \cite{rudin1973fa} for more detailed theory about distributions. We will also use some deeper {results} presented in \cite{knapp2017advanced}.

For the basic terminology, recall that by the support of a tempered distribution $T$, denoted by supp$(T)$, we mean the complement of the largest open set $U$ such that if $\varphi$ is a Schwartz function supported on $U$, {then} $T(\varphi) = 0$. We need the following result, which can be found in \cite[Lemma 4.5]{gabardo2009weighted} and \cite[Theorem 5]{kolountzakis2000structure} (see also \cite[Theorem 7]{kolountzakis2000structure}). 

\begin{lemma}\label{lemma-trans-bdd}
Let $\Lambda$ be translationally-bounded. 
\begin{enumerate}  
\item Suppose that for some $\tau>0$, 
    $$
    \mbox{supp}(\ft\delta_{\Lambda}) \cap B(0,\tau) = \{0\}.
    $$
    Then $\ft\delta_{\Lambda} = a\delta_0$ in $B(0,\tau)$ for some {$a\ge0$}
    \item Suppose that $\ft\delta_{\Lambda}$ is a measure in a neighborhood of $0$. Then $ \ft\delta_{\Lambda}(\{0\}) =\mbox{dens}(\Lambda)$.
\end{enumerate}
\end{lemma}

Around 2010, Gabardo initiated a series of general study on Beurling densities with bounds of the convolutional inequality of the form $\delta_{\Lambda}\ast f$ \cite{Gab12,Gab13}. In particular, among many other results {he proved}, we need the following:

\begin{theorem} \label{Th-Gab} \cite[Proposition 3, Theorem 1 and Corollary 3]{Gab13}
\begin{enumerate}
\item  $D^{+}(\Lambda)<\infty$ if and only if $\Lambda$ is translationally-bounded. 
\item  Let ${0 \le }f\in L^1(\RR^d)$.  Suppose that there exists $C>0$ such that $\delta_{\Lambda}\ast f(x) \le C$ for almost all $x\in \RR^d$. Then 
$$
D^{+}(\Lambda) \cdot \left(\int {f(x)}\,dx\right) \le C. 
$$
\end{enumerate}
\end{theorem}

\subsection{Density properties of tight-frame spectrum}   Suppose $\mu$ is a probability measure on $\RR^d$ which is spectral with $\Lambda \subseteq \RR^d$ being a spectrum. By a well-known result in \cite{jorgensen1998dense}, {this is equivalent to}
\begin{equation}\label{spectrum-as-tiling}
|\ft\mu|^2\ast \delta_{\Lambda} (x) = \sum_{\lambda\in\Lambda}\Abs{\ft{\mu}}^2(x-\lambda) = 1,
\end{equation}
for all $x \in \RR^d$. Note that a spectrum $\Lambda$ must be separated by the continuity of $\ft{\mu}$, $\ft{\mu}(0) = 1$ and  $\widehat{\mu}(\lambda-\lambda') = 0$ for distinct $\lambda,\lambda'\in\Lambda$. This also implies that $\Lambda$ must be translationally bounded.

If $\Abs{\ft\mu}^2\in L^1(\RR^2)$,  then $\mu$ is an absolutely continuous measure with a density in $L^2(\RR^d)$. On the other hand, if the support of $\mu$ has Lebesgue measure 0, then $\int\Abs{\ft\mu}^2 = +\infty$. The following proposition shows that {the density of $\Lambda$} must be zero.

\begin{proposition}\label{prop-Density}
    Suppose that $\mu$ is a singular spectral measure with a {tight-frame} spectrum $\Lambda$. Then the density of $\Lambda$ is zero. i.e. $D^{+}(\Lambda) = 0$. 
\end{proposition}

    \begin{proof}
        In \cite{HLL2013}, it was proved that if $\Lambda$ is a spectrum of $\mu$ we have that $D^-(\Lambda) = 0$. However, to show that it has density zero, we need to  show a stronger statement that $D^{+}(\Lambda) = 0$.\footnote{A sketch of this proof was first shown to the second-named author by J.P Gabardo back in 2012.}

To prove our conclusion. We first notice that $\Abs{\ft\mu}^2*\delta_\Lambda = A$. Let $f_n = \Abs{\ft\mu}^2 \cdot \one_{B_n(0)}$ which is in $L^1$ since $\Abs{\ft\mu} \le 1$. Since $f_n \le \Abs{\ft\mu}^2$ we have $f_n*\delta_\Lambda \le A$. By Theorem \ref{Th-Gab}, we conclude that $D^+(\Lambda) \cdot \int f_n \le A$. Since $\int f_n \to \int \Abs{\ft\mu}^2 = +\infty$, it follows that $D^+(\Lambda)=0$.
 \end{proof}

\begin{proposition}
    Suppose that $\mu$ admits a tight frame of exponentials with tight-frame spectrum $\Lambda$. Then $\Lambda$ is translationally bounded.
\end{proposition}

\begin{proof}
    The previous proposition already showed that $D^{+}(\Lambda) =0$ and hence translationally bounded by Theorem \ref{Th-Gab}(1). Suppose  that $\mu$ is absolutely continuous with respect to the Lebesgue measure, the fact that $\mu$ admits a tight frame of exponentials with tight-frame spectrum $\Lambda$ implies that we have $|\ft\mu|^2\ast \delta_{\Lambda} (x) = A$. By Theorem \ref{Th-Gab} (2), $D^{+}(\Lambda)<\infty$. 
\end{proof}

\subsection{Distributions supported on straight lines.}\label{subsection-line-distribution} Theorem \ref{prop-non-spectral} requires $\mu\ast\wt{\mu}$ to have a smooth density. This is useful if $\mu$ is supported on a smooth curve. However, this is not the case for measures supported on closed polygonal lines. In this subsection, we will prove a lemma about distributions supported on a straight line.

\begin{lemma}\label{lm:no-transverse-derivatives}
Suppose $F \in L^\infty(\RR^2)$ and $T = \ft F$, a tempered distribution, has $\supp T \subseteq \RR\times\Set{0}$.
Then\\
(a) there exists a distribution $\wt T$ on $\RR$ such that for any $h \in \mathcal{S}(\RR^2)$ we have
$$
T(h) = \wt T(h(\cdot, 0)), \text{ and }
$$
(b) $F$ does not depend on $x_2$.
\end{lemma}

\begin{proof}
We can write \cite[Problem 5.6.10]{knapp2017advanced} 
\begin{equation}\label{decomposition}
T(h) = \sum_{j=0}^J T_j(\partial_{0,j}h |_{x_2=0})
\end{equation}
for any Schwartz function $h$ on $\RR^2$. Here $T_j$ are tempered distributions on $\RR$ and $J$ is a finite number. To prove the Lemma it suffices to show that $J=0$.

Let $\phi$ be a smooth compactly supported function of integral 1 and write
$$
\phi_\epsilon(x_1, x_2) =  \phi(x_1, x_2/\epsilon),
$$
which implies
$$
\ft{\phi_\epsilon}(\xi_1, \xi_2) = \epsilon \ft{\phi}(\xi_1, \epsilon \xi_2).
$$
We have
\begin{equation}\label{bounded}
\Abs{T(\phi_\epsilon)} = \Abs{F(\ft{\phi_\epsilon})} = \Abs{ \int_{\RR^2} F \ft{\phi_\epsilon} } < \Linf{F} \Abs{\phi(0)},
\end{equation}
independent of $\epsilon$. 

Using \eqref{decomposition} and the fact that
$$
\partial_{0,j}\left(\phi_\epsilon(x_1, x_2)\right) =
 \epsilon^{-j} (\partial_{0,j}\phi) (x_1, x_2/\epsilon)
$$
we get
$$
T(\phi_\epsilon) = \sum_{j=0}^J \epsilon^{-j} T_j((\partial_{0,j} \phi)(x_1, 0)).
$$
Assume now that $J>0$ and $\partial_J T_J \neq 0$ (otherwise the term does not exist) and choose a Schwartz function $\phi$ for which $T_J((\partial_{0,j} \phi)(x_1, 0)) \neq 0$. It follows that as $\epsilon \to 0+$ the quantity $\Abs{T(\phi_\epsilon)}$ grows as $\epsilon^{-J}$, which contradicts \eqref{bounded}, and concludes the proof of (a).

To prove (b) it is enough to show that $\tau_{(0,\delta)}F = F$, where the translation operator $\tau_h$ acts on a distribution $S$ as follows:
$$
(\tau_{(h_1, h_2)} S)(\phi(x_1, x_2)) = S(\phi(x_1-h_1, x_2-h_2).
$$
From part (a) we have 
\begin{align*}
\ft{\tau_{(0, \delta)} F} (h) &= (e^{2\pi i \delta \xi_2} \ft{F} )(h) \ \text{ (translation becomes modulation) }\\
& = \ft{F} ( e^{2\pi i \delta \xi_2} h(\xi_1, \xi_2)) \\
& = \wt T(h(\xi_1, 0)) \\
&= \ft{F} (h).
\end{align*}
By Fourier uniqueness we thus conclude that $\tau_{(0, \delta)} F = F$.
\end{proof}

\subsection{A change of variable formula.}     Let $F: \RR^d\to\RR^d$ be a Lipschitz function. The Jacobian of $F$ is defined almost everywhere as usual and will be denoted by $J_F$.  We need the following change of variable formula in geometric measure theory, which is known as the {\it area formula} (see \cite[{\S 3.3.3}]{Evans-Gariepy}). 

\begin{theorem}\label{area formula}
    Let $F: \RR^d\to\RR^d$ be a Lipschitz function. Then for all $g\in L^1(\RR^d)$. 
    $$
    \int g(x)|J_F(x)|~dx =  \int \left(\sum_{x\in F^{-1}\{u\}} g(x)\right) ~du
    $$
\end{theorem}
In particular, If $F$ is a bijective function, we have the change-of-variable formula of the following forms:
\begin{equation}\label{eq-change-1}
    \int g(F(x))|J_F(x)|~dx = \int g(u)~du
\end{equation}
and 
\begin{equation}\label{eqchange2}
    \int g(F(x))~dx = \int  g(u) | J_{F}(F^{-1}(u))|^{-1}~du.
\end{equation}
(\ref{eq-change-1}) is our usual change of variable formula, which can be seen from the area formula. (\ref{eqchange2}) follows from (\ref{eq-change-1}) by noting that $|J_F|^{-1}|J_F| = 1$. 


\section{Proofs of Theorems \ref{prop-non-spectral} and \ref{prop-spectral-gap}} \label{s:thm34} We will give the proof for Theorem \ref{prop-non-spectral} and \ref{prop-spectral-gap} in the introduction. Combined with the preliminary tools, they form the basic framework to the proof of Theorem \ref{theorem1} and \ref{theorem2}

\begin{proof}  [Proof of Theorem \ref{prop-non-spectral}] Let $\psi$ be a smooth function of compact support in $U$ and define $h = \psi\cdot (\mu*\wt\mu)$, which is smooth everywhere by our assumption. Then $\tf h = \tf\psi * \Abs{\ft\mu}^2$ so
$$
\tf h*\delta_\Lambda = \tf\psi*\Abs{\ft\mu}^2*\delta_\Lambda = \tf\psi * 1 = \psi(0).
$$
It follows that $h \cdot \ft\delta_\Lambda = \psi(0)\delta_0$ which implies that
$$
\supp\ft\delta_\Lambda \cap \Set{h \neq 0} \subseteq \Set{0}.
$$
Since $\Set{h \neq 0} = \Set{\psi \neq 0}$ and the latter set can be taken to be an open neighborhood of any point in $U$ the desired inclusion follows.
\end{proof}

 \begin{proof} [Proof of Theorem \ref{prop-spectral-gap}].
     We argue by contradiction. Suppose that $\Lambda$ is a {tight-frame} spectrum for some singular measure $\mu$ and $\delta_{\Lambda}$ has a spectral gap. Then there exists an open set $U$ such that \begin{equation}\label{support-at-0}
\supp\ft\delta_\Lambda \cap U \subseteq \Set{0}.
\end{equation}
     Note that a spectrum $\Lambda$ must be translationally bounded. By Lemma \ref{lemma-trans-bdd} (1) and equation (\ref{support-at-0}), we have $\ft\delta_{\Lambda} = a\delta_0$ in $B(0,\tau)$. This shows that $\ft\delta_{\Lambda} $ is locally a measure. Lemma \ref{lemma-trans-bdd} (2) implies that $a = \dens(\Lambda)$ and this is equal to zero by Proposition \ref{prop-Density}. Hence, $\ft\delta_{\Lambda} $ is locally zero around the origin. 

\begin{figure}[h]
\ifdefined\SMART\resizebox{4cm}{!}{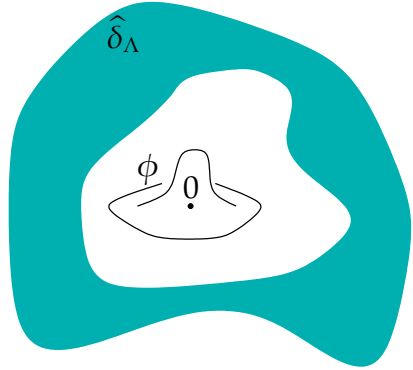}\else
\input test-function.pdf_t
\fi
\caption{A test function $\phi$ near the origin with $\int\phi>0$ leads to a contradiction in the proof of Theorem \ref{prop-spectral-gap}.}\label{fig:test-function}
\end{figure}

     However, this contradicts the fact that $\delta_\Lambda$ is a nonnegative measure. Indeed take a smooth function $\phi$, with $\int\phi>0$, compactly supported in a sufficiently small neighborhood of the origin (see Fig.\ \ref{fig:test-function}). Then we  have $\ft\delta_\Lambda(\phi)=0$. On the other hand, $\ft\delta_\Lambda(\phi) = \delta_\Lambda(\ft\phi) = \sum_{\lambda\in\Lambda}\ft\phi(\lambda) \ge \ft\phi(0) > 0$, a contradiction.
\end{proof}

     We also notice that Theorem \ref{prop-spectral-gap} can also be deduced from a result in \cite[Proposition 4]{LO15}.

\section{Non-spectrality of {collections of} line segments}\label{s:line}
\subsection{Geometric setup.} We need to set up some notation for our discussions. Let ${\bf v}$ and ${\bf w}$ be two linearly independent vectors in $\RR^2$. The parallelogram generated by ${\bf v}$ and ${\bf w}$ will be denoted by 
$$
Q_{\bf v,w} = \{t{\bf v}+s{\bf w}: t,s\in{[0,1]}\}.
$$
For a given ${\bf x}_0\in \RR^2$, the line segment beginning with ${\bf x}_0$ in the direction of ${\bf v}$ is given by 
$$
L_{{\bf x}_0,{\bf v}} = \{x_0+t{\bf v}: t\in[0,1]\}.
$$
We will endow $L_{{\bf x}_0,{\bf v}}$ with the 1-dimensional Hausdorff measure $\mu$ given by 
\begin{equation}\label{eq-line-measure}
\int f(x)~d\mu (x) = \int_0^1 f({\bf x}_0+t{\bf v})~dt
\end{equation}
for all integrable functions $f.$
\begin{lemma}\label{lemma_convolution}
Let $\mu,\nu$ be the line measure supported respectively on the line ${\bf x}_0+t{\bf v}$ and ${\bf y}_0+t{\bf w}$, with $t\in[0,1]$, given by (\ref{eq-line-measure}).
    \begin{enumerate}
    \item Suppose that ${\bf v}$ and ${\bf w}$ are not parallel. Then $\mu\ast\wt{\nu}$ is, up to a constant, the 2-dimensional Lebesgue measure supported on  ${\bf x}_0-{\bf y}_0 + Q_{\bf v,-w}$. 
    \item Suppose that ${\bf v}$ and ${\bf w}$ are parallel and write ${\bf w} = a{\bf v}$. Then $\mu\ast\wt{\nu}$ is supported on the line $\{({\bf x}_0-{\bf y}_0)+ t{\bf v}: t\in\RR\}$ with a continuous density with respect to the one-dimensional Hausdorff measure on the line.
    \end{enumerate}
\end{lemma}

\begin{proof}
\noindent (1). By the definition of convolution and (\ref{eq-line-measure}), 
$$
\int f d{(\mu\ast\wt\nu)} = \iint f(x-y) d\mu(x)d\nu (y) = \int_0^1\int_0^1 f({\bf x}_0-{\bf y}_{0}+ t{\bf v}-s{\bf w}) ~dtds.
$$
Applying the change of variable formula (\ref{eqchange2}) we see that the above is equal to 
$$
c^{-1}\int f(u) du, \ \mbox{where} \  c = \det\left({\bf v}, -{\bf w}\right)\ne 0.
$$
By checking $f = \one_E$, we see that $\mu\ast\wt\nu(E)>0$  if and only if $E\subset {\bf x}_0-{\bf y}_0 + Q_{\bf v,-w}$, and it is equal to $c^{-1}$ times the Lebesgue measure of $E$ whenever $E$ is a subset of the parallelogram. This shows (1).

\noindent (2). In a similar calculation, if ${\bf w}  = a{\bf v}$, 
$$
\int f d{(\mu\ast\wt\nu)} = \iint f(x-y) d\mu(x)d\nu (y) = \int_0^1\int_0^1 f({\bf x}_0-{\bf y}_{0}+ (t-as){\bf v}) ~dtds.
$$
Checking with characteristic functions, the measure is supported on the line $L_{{\bf x}_0-{\bf y}_{0}-a{\bf v},(1+a){\bf v}}$. We now compute its density.  Letting $g(x) = f({\bf x}_0-{\bf y}_{0}+x{\bf v})$ and $m_a$ to be the Lebesgue measure on $[0,a]$, we see that the above integral can be rewritten as
$$
\int_0^1\int_0^1 g(t-as)~dtds = \int g ~dm_1\ast m_a.
$$
As $m_1\ast m_a$ has a continuous density, the conclusion of (2) follows. 
\end{proof}

In the following, a closed half-plane $H$ in $\RR^2$ is the set $\{x\in\RR^2: \langle{\bf u}, {\bf x}\rangle\ge 0\}$ for some vector ${\bf u} \neq {\bf 0}$. 
\begin{lemma}\label{lemma-covering}
    Let $N\ge 3$ and let  ${\bf v}_1,\cdots,{\bf v}_N$ be non-zero vectors such that there does not exist a closed half-plane $H$ such that ${\bf v}_i\in H$ for all $1\le i\le N$. Then 
    $$
    \bigcup_{1\le i\ne j\le N}Q_{{\bf v}_i,{\bf v}_j}
    $$
    covers an open set containing the origin. 
\end{lemma}

\begin{proof}
We will prove it by induction on $N$. First, we prove the case when $N=3$. Let ${\bf v}_1,{\bf v}_2,{\bf v}_3$ be three non-zero vectors. By our assumption, none of the two vectors can be parallel since otherwise, it will be contained in the same half-plane. Taking ${\bf v}_1,{\bf v}_2$. Then they are linearly independent and 
$$
{\bf v}_3 \in \{t{\bf v}_1+s{\bf v}_2: t,s<0\}.
$$
Otherwise, ${\bf v}_1,{\bf v}_2,{\bf v}_3$ will lie on the same half-plane. Now, it is a routine check that the parallelograms these three vectors generated cover a neighborhood of origin. 

\medskip

Suppose that the statement is true for $N-1$. We now prove it also holds for $N$. Indeed, given non-zero vectors  ${\bf v}_1,\cdots,{\bf v}_N$. Suppose that the first $N-1$ vectors does not lie on any half plane. Then we can apply induction hypothesis to obtain our desired conclusion. Thus, we assume that the first $N-1$ vectors lies in some half-plane. By a rotation, we can assume that the half-plane is $y\ge 0$ and assume that the vectors ${\bf v}_1,\cdots,{\bf v}_{N-1}$ are arranged in anticlockwise directions from the positive $x$-axis.

{\bf Case (1).} Suppose that the angle between ${\bf v}_1$ and ${\bf v}_{N-1}$ is strictly less than $\pi$. Then all ${\bf v}_{i}$,  $i=2,\cdots N-2$ will lies in the quadrant $\{t{\bf v}_1+s{\bf v}_2: t,s\ge 0\}$. Hence, 
$$
{\bf v}_N \in \{t{\bf v}_1+s{\bf v}_{N-1}: t,s<0\}.
$$
Otherwise, all vectors will lie in the same half plane. Then we consider the three vectors ${\bf v}_1,{\bf v}_{N-1},{\bf v}_N$ and apply the case for three vectors, we conclude that $Q_{{\bf v}_1,{\bf v}_{N-1}}\cup Q_{{\bf v}_{N-1},{\bf v}_{N}}\cup  Q_{{\bf v}_{N},{\bf v}_{N}}$ covers a neighborhood of the origin. 

{\bf Case (2).} Suppose that the angle between ${\bf v}_1$ and ${\bf v}_{N-1}$ is  $\pi$. In this case, ${\bf v}_N$ must belong to $y<0$. We can also find ${\bf v}_j$ for some $1<j<N-1$ such that ${\bf v}_j$ lies in the plane $y>0$. Hence, the parallelograms generated by ${\bf v}_1,{\bf v}_j,{\bf v}_{N-1}, {\bf v}_N$ provides a covering of the neighborhood of the origin. This completes the whole proof.

\end{proof}

{Alternatively, we notice that the assumption about the half-plane is equivalent to 0 being in the convex hull of the points. After finishing the proof of three vectors, by the {Caratheodory's theorem} in convex geometry, 0 must belong to the convex hull of three of the vectors, and since the union in this lemma contains the corresponding union for the 3 vectors, we conclude the desired result.}

\subsection{Closed polygonal lines.} Let ${\bf v}_1,\cdots, {\bf v}_N$ be $N$ non-zero vectors in $\RR^2$ such that 
$$
\sum_{i=1}^N {\bf v}_i = {\bf 0}.
$$
Since we can combine two adjacent parallel vectors together, we may assume without loss of generality that {successive ${\bf v}_i$} are not parallel.  Let ${\bf x}_0 = {\bf 0}$, ${\bf x}_i = {\bf v}_1+\cdots+{\bf v}_i$. Define  the line 
$$
L_i = \{{\bf x}_{i-1}+t{\bf v}_i: t\in[0,1]\}, \ i = 1,\cdots, N.
$$ 
These {line segments} form a closed polygonal curve ({which} may cross itself). Endow each $L_i$ with Lebesgue measure $\mu_i$ so that 
$$
\int f(x)d\mu_i(x) = \int_0^1 f({\bf x}_{i-1}+t{\bf v}_i) ~dt.
$$
Define the measure on the closed polygonal line
$$
\mu = \frac1N \sum_{i=1}^N \mu_i.
$$ 
Let  ${\mathcal L}_i: = \{-{\bf v}_i t: t\ge 0\}$ be the half-line. Note that the vectors ${\bf v}_i$ and ${\bf v}_j$ may be in the same parallel vectors.  ${\mathcal L}_i ={\mathcal L}_j$ is the same if and only if  ${\bf v}_i = a {\bf v}_j$ for some $a>0$. Let $M$ be the number of all these distinct half-lines. Define also 
$$
{\mathcal W}_i : = \{t{\bf v}_i: t\in\RR\}.
$$
\begin{lemma} \label{lemma-structure-mu*mu} Let  $\mu$ be the measure defined above. Then there exists a bounded open neighborhood of the origin $U$ such that  $U$ is divided into $M$ regions by the {half-}lines ${\mathcal L}_i, i = 1,\cdots,M$, and  
$$
 \mu\ast\wt{\mu} = {F(x)}+ \nu_1+\cdots+\nu_M
 $$
 where
 \begin{enumerate}
 \item $\nu_i$ is a singular measure compactly supported on the subspace ${\mathcal W}_i$ with $0$ is in the support; and 
 \item $F$ is a piecewise constant function, which is a constant in each of the regions. 
 \end{enumerate}
 \end{lemma}

\begin{proof}
  First we know that 
  $$
 \mu\ast\wt{\mu} = \frac1{N^2}\sum_{i=1}^N\mu_i\ast\wt\mu_i+ \frac{1}{N^2}\sum_{1\le i\ne j\le N} \mu_i\ast\wt\mu_j.
$$
By Lemma \ref{lemma_convolution}(2), $\mu_{i}\ast\wt{\mu}_i$ is a singular measure supported on the subspace ${\mathcal W}_i$ with $0$ is in the support. This gives the measure $\nu_i$. 

We look at $\mu_i\ast\wt\mu_j$ for $i\ne j$. In particular, let us consider $j = i+1$ in which the vectors ${\bf v}_i,{\bf v}_{i+1}$ are not parallel. By Lemma \ref{lemma_convolution} (1), $\mu_i\ast\wt\mu_{i+1}$ is a constant function on 
$$
{\bf x}_{i-1}-{\bf x}_i+ Q_{{\bf v}_i,-{\bf v}_{i+1}} = Q_{-{\bf v}_i,-{\bf v}_{i+1}},
$$
where the last set equality follows from the definition of ${\bf x}_i$. We notice that as the sum of $-{\bf v}_i$ is a zero vector, $\{-{\bf v}_i\}_{i=1}^N$ cannot lie on the same half-plane. By Lemma \ref{lemma-covering}, the union of $Q_{-{\bf v}_i,-{\bf v}_j}$ covers a neighborhood of the origin, denoted by $U$. As in each region determined by the lines ${\mathcal L}_i$, it is covered by finitely many $Q_{-{\bf v}_i,-{\bf v}_j}$. It must therefore be a constant function. 
\end{proof}

\begin{proof}[Proof of Theorem \ref{theorem1}(1)]\label{ss:support} We are now ready to show that $\mu$ cannot be a tight-frame spectral measure.  We consider the open neighborhood $U$ in Lemma \ref{lemma-structure-mu*mu}. Away from the subspaces ${\mathcal W}_i$, $\mu\ast\wt\mu$ is a constant function. By Theorem \ref{prop-non-spectral}, 
$$
\supp\ft\delta_\Lambda \cap U \subseteq \bigcup_{i=1}^N{\mathcal W}_i.
$$
\begin{center}
    {\bf Claim:} $\ft\delta_\Lambda$ has no support on all ${\mathcal W}_i$ except possibly at the origin.
\end{center}

To see this claim, we fix a subspace ${\mathcal W}_i$ and by a rotation, we may assume that ${\mathcal W}_i$ is the $x$-axis.  By taking the function $f \in L^2(\mu)$ so that $f(x,0) = e^{2\pi t x}$ and be equal to $0$ on the other sides, and applying the  Parseval's identity,  we obtain that
$$
\sum_{\lambda = (\lambda_1, \lambda_2) \in \Lambda} \Abs{\ft{\one_{[0, 1]}}}^2(t-\lambda_1) = N.
$$
This tiling condition implies that $\Lambda_1 = \pi_1\Lambda$ (the projection of $\Lambda$ onto the first axis, a multiset) has density $N$.

Let $\epsilon>0$ be sufficiently small and take a smooth $\phi$ supported in an $\epsilon$-disk around the origin.
Define $\psi(x_1, x_2) = \phi((x_1, x_2)-(\epsilon, 0))$. The claim will follow if we can show that $\psi \ft{\delta_\Lambda}$ is the zero distribution. Suppose it is not, so that its inverse Fourier transform $F = \tf\psi*\delta_\Lambda$ is a non-zero bounded function on $\RR^2$ since $\tf\psi$ is a Schwartz function and $\Lambda$ is translationally bounded (see Lemma \ref{lemma-schwartz-bounded}).

Since $F$ is a \textit{complex} function that we assume not to be identically 0 there is no loss of generality to assume that
$\Re F(x_1, 0) \ge c > 0$ for $x_1 \in E \subseteq (-M, M)$, some bounded measurable subset of $\RR$ of positive measure. By Lemma \ref{lm:no-transverse-derivatives} (b), $F$ is constant in $x_2$.
 Then $\int_{E \times \RR} \Re F = +\infty$. We have 
$$
+\infty = \int_{E \times \RR} \Re F = \int_{E \times \RR} \Abs{\Re F} \leq \int_{E \times \RR} \Abs{F}.
$$
Expanding out the definition of $F$, this implies that 
$$
\int_{E \times \RR} \Abs{F} \le \int\int_{(-M, M) \times \RR} \sum_{(\lambda_1,\lambda_2)\in \Lambda} |\tf\psi(x_1-\lambda_1,x_2-\lambda_2)|dx_1dx_2.
$$
Recall that $\Lambda_1= \pi_1(\Lambda)$ and we define $h(x_1) = \int_{\RR} |\tf\psi(x_1,y)|dy$. Then using the translatioal invariance of Lebesgue measure in the $x_2$ coordinates, 
$$
\int\int_{(-M, M) \times \RR} \sum_{(\lambda_1,\lambda_2)\in \Lambda} |\tf\psi(x_1-\lambda_1,x_2-\lambda_2)|dx_1dx_2 = \int_{(-M,M)} \sum_{\lambda_1\in\Lambda_1} |h(x_1-\lambda_1)| dx_1
$$
Finally, notice that $h$ is a Schwartz function and $\Lambda_1$ has a translationally bounded, so $\sum_{\lambda_1\in\Lambda_1} |h(x_1-\lambda_1)|$ is a bounded function with some bound $C<\infty$ (by Lemma \ref{lemma-schwartz-bounded}).  This shows that the above integral is at most $2CM$. This is a contradiction and this justifies the claim. 

With the claim justified, we can conclude that $\supp\ft\delta_\Lambda \cap U \subseteq \{0\}$. Hence, $\delta_{\Lambda}$ has a spectral gap. This contradicts Theorem \ref{prop-spectral-gap}.  Hence, $\mu$ cannot be a spectral measure and this completes the proof.
\end{proof}

\subsection{Line segments from the same points.} We will prove Theorem \ref{theorem1} (2) in this subsection. Fix ${\bf x}_0\in \RR^2$ and {let} ${\bf v_1},\cdots, {\bf v}_N$ be non-zero vectors. We consider the {line segments}
$$
L_i = \{{\bf x}_0+t{\bf v}_i: t\in[0,1]\}
$$
and $\mu_i$ {is} the Lebesgue measure on $L_i$. Theorem \ref{theorem1} (2) requires us to show that the measure $\mu = \frac{1}{N} \sum_{i=1}^N \mu_i$ is not a spectral measure. The proof follows from the same strategy {as} in Theorem \ref{theorem1} (1).

\begin{lemma}\label{lemma-covering2}
    Let  ${\bf v}_1,{\bf v}_2,{\bf v}_3$ be non-zero vectors such that either
    \begin{enumerate}
        \item ${\bf v}_3 = a{\bf v}_1$ for some $a<0$ and ${\bf v}_1$ and ${\bf v}_2$ are not parallel or 
         \item   they are not pairwise parallel,
    \end{enumerate}
 Then 
    $$
    \bigcup_{1\le i\ne j\le 3}Q_{{\bf v}_i,-{\bf v}_j}
    $$
    covers an open set containing the origin and  the interior of all $Q_{{\bf v}_i,-{\bf v}_j}$ are pairwise disjoint. 
\end{lemma}

\begin{proof} We use ${\bf v}_1$ and ${\bf v}_2$ to span $\RR^2$. It divides $\RR^2$ into 4 quadrants. We notice that $Q_{{\bf v}_1,-{\bf v}_2}$ and $Q_{{\bf v}_2,-{\bf v}_1}$ already cover the first and the third quadrant around the origin. If ${\bf v}_3 = a{\bf v}_1$ for some $a<0$, then  $Q_{{\bf v}_2,-{\bf v}_3}$ and $Q_{{\bf v}_3,-{\bf v}_2}$ cover the second and the fourth quadrants. Hence, their union covers a neighborhood of the origin and they are pairwise disjoint.  This proves the first case.

In the second case, we subdivide into four cases according to  ${\bf v}_3$ lying strictly inside each quadrant. By a standard check, the conclusion also holds. 
\end{proof}

\begin{proof}[Proof of Theorem \ref{theorem1} (2)]  Suppose that the union of lines contains three non-parallel lines from ${\bf x}_0$ or two parallel vectors in opposite directions and the last one non-parallel. Suppose that this measure is a spectral measure with a spectrum $\Lambda\subset \RR^2$. By restricting to the three lines generated by ${\bf v}_1,{\bf v}_2,{\bf v}_3$, we consider $\rho = \frac13(\mu_1+\mu_2+\mu_3)$ where $\mu_i$ is the measure supported on the {line segment} $\{{\bf x}_0+t{\bf v}_i: t\in[0,1]\}$. As $\Lambda$ is a spectrum for $\mu$, considering only functions supported on these three lines,
$$
\sum_{\lambda\in\Lambda} \left|\int f(x)e^{-2\pi i \lambda\cdot x}d\rho (x)\right|^2  = \frac{N^2}{3^2} \|f\|^2_{L^2(\rho)}.
$$
By taking $f = e^{2\pi i \xi\cdot x}$, we have $\delta_{\Lambda}\ast |\ft{\rho}|^2 = w$, where $w= N^2/9$. We now compute 
  $$
 \rho\ast\wt{\rho} = \frac1{3}\sum_{i=1}^N\rho_i\ast\wt\rho_i+ \frac{1}{9}\sum_{1\le i\ne j\le 3} \rho_i\ast\wt\rho_j.
$$
By Lemma \ref{lemma_convolution} (2), $\rho_{i}\ast\wt{\rho}_i$ is a singular measure supported on the subspace ${\mathcal W}_i = \{t{\bf v}_i: t\in\RR\}$ with $0$ in the support. Using Lemma \ref{lemma-covering2}, the support of $\sum_{1\le i\ne j\le 3} \rho_i\ast\wt\rho_j$ is exactly 
$$
\bigcup_{1\le i\ne j\le 3} Q_{{\bf v}_i,-{\bf v}_j}, 
$$
and the measure $\rho_i\ast\wt\rho_j$ is absolutely continuous with  a constant density on each  $Q_{{\bf v}_i,-{\bf v}_j}$. Hence,  Using Theorem \ref{prop-non-spectral}, we conclude that in an open set $U$ containing $0$
$$
\supp \ft{\delta_{\Lambda}}\cap U \subset \{0\}\cup {\mathcal W}_1\cup{\mathcal W}_2\cup {\mathcal W}_3.
$$
We see that the claim in the proof of Theorem \ref{theorem1} (1) in the previous subsection works analogously here.  We conclude that $\supp \ft{\delta_{\Lambda}}\cap U \subset \{0\}$. This leads to a contradiction to Theorem \ref{prop-spectral-gap} as we now have a spectral gap for $\delta_{\Lambda}$.
\end{proof}


\subsection{Two particular cases of shapes and their $\mu*\wt \mu$ measure.}\label{ss:two-cases}
\begin{figure}[h]
\ifdefined\SMART\resizebox{4cm}{!}{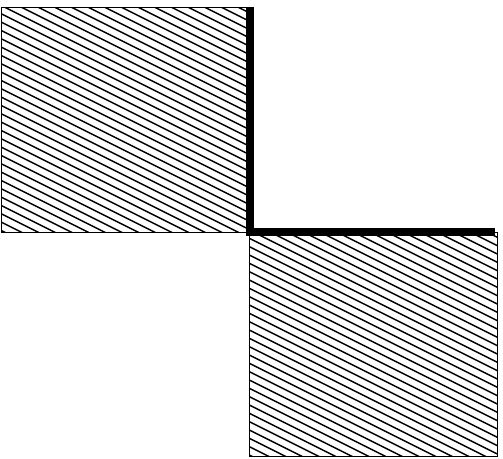}\else
\input l-shape.pdf_t
\fi
\caption{The ``L-shape'', the union of the two line segments from the origin to $(1, 0)$ and to $(0, 1)$, is spectral \cite{lai2021spectral}.}\label{fig:l-shape}
\end{figure}

It is also instructive to exhibit one case where our technique does not apply and the measure is indeed spectral. It is the ``L-shape'', shown in Fig.\ \ref{fig:l-shape}.

 The measure $\mu*\wt \mu$ when $\mu$ is arc-length on the L-shape, is supported on the shaded area in Fig.\ \ref{fig:l-shape} and it also has a singular part supported on the two line segments from $(-1, 0)$ to $(1, 0)$ and from $(0, -1)$ to $(0, 1$). Thus $\mu*\wt \mu$ does not cover a neighborhood of the origin and the proof of Theorem \ref{theorem1} does not apply.

 Having seen the measure $\mu*\wt \mu$ for the L-shape we can easily see from this that the corresponding convolution measure for the $\Pi$-shape of Fig.\ \ref{fig:polygons} (the polygonal line joining $(0, 0)$ to $(0, 1)$ to $(1, 1)$ to $(1, 0)$) has a smooth part supported on the open set that we get if we remove the coordinate axes from $(-1, 1)^2$ plus a singular part supported on the same two line segments as the L-shape, namely the line segments from $(-1, 0)$ to $(1, 0)$ and from $(0, -1)$ to $(0, 1$). It follows from the same proof as that of Theorem \ref{theorem1} (1) that the $\Pi$-shape is not spectral.

\section{Non-spectrality of Curves}\label{Section-curves}

A smooth closed curve is an infinitely differentiable function $\gamma: [a,b]\to \RR^2$ such that all derivatives $\gamma^{(n)}(a) = \gamma^{(n)}(b)$ for all $n\ge 0$.  We can always parametrize $\gamma$ by arc length so that  $\gamma:[0,L]\to \RR^2$ be a smooth curve with the arc length parametrization ($L$ is the arc length of $\gamma$). Let ${\bf n}$ be the unit normal vector so that ${\bf n}$ is orthogonal of $\gamma'$ and $\{\gamma',{\bf n}\}$ is of positive orientation. Note that $\gamma''$ is parallel to ${\bf n}$.  The curvature of $\gamma$ is defined to be 
$$
\kappa(s) = \pm \|\gamma''(s)\|
$$
where the sign  is determined by whether $\gamma''(s)$ is in same or opposite direction of ${\bf n}$. We say that a curve has positive curvature if its curvature never vanishes. A consequence of this is that, as the parameter $s$ increases, the tangent vector $\gamma'(s)$ is always turning in the same direction, either clockwise or counterclockwise. For more detailed theory about the geometry of curves, readers are invited to consult \cite{docarmo}.  

In our study of curves, we do not assume the curve is simple {(i.e. not self-intersecting)}. Suppose that $x_0$ is an intersection point of the curve with itself. We say that the intersection point $x_0$ is {\bf transverse} if there exists only finitely many preimages of $x_0$ under $\gamma$ and  $\gamma' (s)\ne \gamma'(t)$ for all distinct $s,t$ such that $\gamma(s) = \gamma(t) = x_0$.

The {\bf arc-length measure} of $\gamma$ is the measure such that 
$$
\int f(x)~d\mu(x) = \int_0^L f(\gamma(t))~dt
$$
for all continuous functions $f$ in $\RR^2$. Without loss of generality, by rescaling the curve, we may assume that $L = 1$, so that $\mu$ is a probability measure.  
\begin{equation}\label{eq-auto-correlation}
\int f(x)d(\mu\ast\wt{\mu})(x) = \int_0^1\int_0^1 f(\gamma(t)-\gamma(s))~dtds
\end{equation}
{Let the function $\Gamma:[0, 1]^2 \to \RR^2$ be defined by}
\begin{equation}\label{eq-Gamma(s,t)}
\Gamma (s,t) = \gamma(s)-\gamma(t). 
\end{equation}
Suppose that $\gamma(s) = (x(s),y(s))$. Then the Jacobian of $\Gamma$ can be computed easily as 
\begin{equation}\label{eq_Jacobian}
J_{\Gamma}(s,t) = \det \begin{pmatrix} x'(s) & -x'(t)\\y'(s) & -y'(t)\\\end{pmatrix} = x'(t)y'(s)-x'(s)y'(t). 
\end{equation}
The determinant is non-zero if and only if $\gamma'(s)$ is not parallel to $\gamma'(t)$.

\begin{lemma}\label{lemma_self}
    Let $\gamma$ be a curve with positive curvature. There exists $\delta_0>0$ such that for all intervals $I$ of length less than $\delta_0$, $\Gamma$ is injective on $(I\times I)\setminus\{(s,s):s\in I\}$.
\end{lemma}

\begin{proof}
Let us partition $[0,1]$ into intervals of length $N^{-1}$ for some large $N$ so that
{ the angle formed by $\gamma'(s)$ and $\gamma'(s')$ is less than $\frac{\pi}{8}$}
for all $s,s'$ inside an interval of length $N^{-1}$. Then for all intervals {$I = (a, b)$} of length $\delta_0$ less than $(2N)^{-1}$, $\theta(s)$ is strictly increasing over $I$ with an angle at most $\pi/4$, {where $\theta(s)$ denotes the angle formed by $\gamma'(s)$ with $\gamma'(a)$.}

By the generalized mean value theorem and the fact that $\theta(s)$ strictly increasing over $I$, for all $(s,t)\in I\times I$ and $s\ne t$, there exists a unique $\xi\in I$, such that $\gamma(s)-\gamma(t)$ is in the same direction as $\gamma'(\xi)$. This shows that $\Gamma$ must be injective.
\end{proof}

\begin{lemma}\label{lemma-tranverse}
    Let $\gamma$ be a smooth closed curve with positive curvature and let $x_0 = \gamma(s_0) = \gamma(t_0)$ be a transverse intersection point. Then there exist intervals $I_{s_0,t_0} = (s_0-\delta_{s_0,t_0},s_0+\delta_{s_0,t_0})\ni s_0$ and $J_{s_0,t_0}=(t_0-\delta_{s_0,t_0},t_0+\delta_{s_0,t_0})\ni t_0$ such that $\Gamma$ is injective on $I\times J$ and the image $\Gamma(I_{s_0,t_0}\times J_{s_0,t_0})$ covers $B(0, \varepsilon_{s_0,t_0})$ for some $\varepsilon_{s_0,t_0}>0$. 
\end{lemma}

\begin{proof}
    As the intersection is transverse, $\det J_{\Gamma}(s_0,t_0)\ne 0$. It follows that $\Gamma$ is locally bijective around $(s_0,t_0)$, which is exactly the lemma required.
\end{proof}

\begin{lemma}\label{lemma-ball cover}
    Let $\gamma$ be a smooth closed curve with positive curvature. Then there exists $\varepsilon_1>0$ such that $\Gamma$ maps  onto $B(0,\varepsilon_1)$. 
\end{lemma}

\begin{proof}
Assume the curve $\gamma:[0,1] \to \RR^2$ parametrized by arc-length.
Because of the non-zero curvature assumption the unit tangent vector $\gamma'(t)$ of $\gamma$ is turning strictly monotonically counterclockwise when we traverse the curve in one of its two orientations, say when $t$ goes from {$0$ to $1$}. Thus $\gamma'(t)$ takes all possible orientations at least once. {We also remark that the curvature is bounded by compactness.}

\begin{figure}[h]
\ifdefined\SMART\resizebox{4cm}{!}{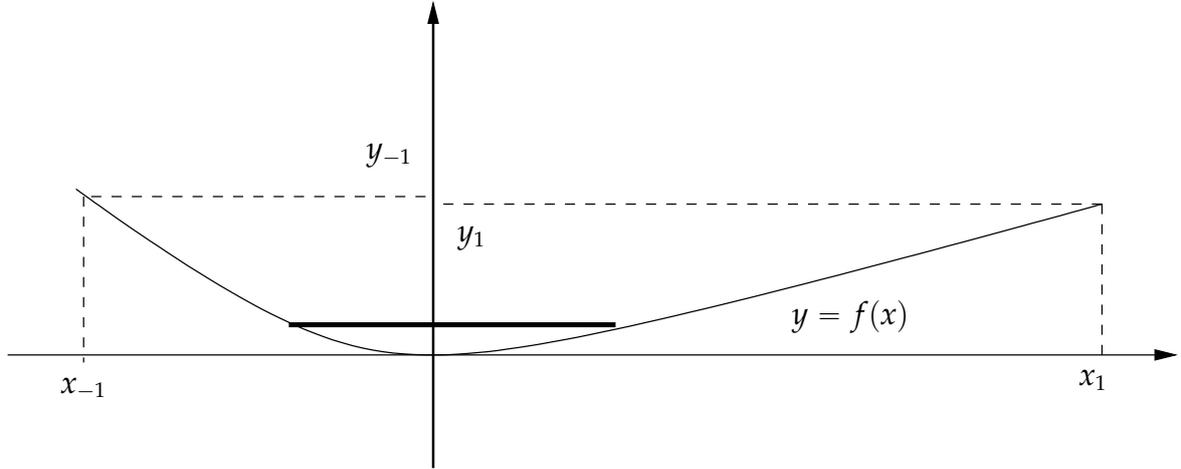}\else
\input horizontal-chord.pdf_t
\fi
\caption{A chord of the right length parallel to a given tangent.}\label{fig:chord}
\end{figure}

Fix a point $p = \gamma(t_0)$ with tangent line {$T$}. We may assume that $p=(0, 0)$ and that the line {$T$} is the $x$-axis. Let $t_{-1} < t_0 < t_1$ be the first points to the left and right of $t_0$ where the tangent has slope $\pm1$. Between $t_{-1}$ and $t_1$ the curve $\gamma$ is the graph of a function $f:[x_{-1}, x_1] \to \RR^+$, where $p(t_{-1}) = (x_{-1}, y_{-1})$ and $p(t_1) = (x_1, y_1)$.
The upper and lower bounds on the curvature of $\gamma$, given by
$$
\kappa(x) = \frac{\Abs{f''(x)}}{\left(1+\Abs{f'(x)}^2\right)^{3/2}},
$$
translate to positive upper and lower bounds on $f''$, independent of $t_0$ (since $\Abs{f'}\le 1$). By the upper bound on $f''$ we obtain a lower bound $\ell$ for $\Abs{x_{-1}}, x_1$ independent of $t_0$. By the lower bound on $f''$ we obtain a lower bound of the form $\Abs{f'(x)} \ge c \Abs{x}$, with $c$ independent of $t_0$ and this in turn gives positive lower bound $h$ on $y_{-1}$ and $y_1$ independent of $t_0$. By looking at the chord of the curve defined by the horizontal line at height $y$ we obtain that all horizontal chords of length from 0 to $2\ell$ (which is bounded below independent of $t_0$) are realized.
This implies that chords of every orientation and length up to $\epsilon_1$ are achievable on the curve $\gamma$.
\end{proof}

\begin{theorem}\label{theorem-smooth}
    Let $\gamma: [0,1]\to \RR^2$ be a closed smooth curve of positive curvature and {assume that} $\gamma$ has {finitely many self-intersections, all of them transverse}. Let $\mu$ be the arc-length measure on $\gamma$. Then $\mu\ast\wt\mu$ is absolutely continuous with a smooth density in a {punctured} neighborhood of $0$ in $\RR^2$. 
\end{theorem}

\begin{proof}
We first note that {since the} curve is closed and all derivatives agree on the end-points, we may assume that the domain of the interval is the circle ${\mathbb T}$ where $0$ and $1$ are identified as the same point, so that it is a compact set.

    Let $\delta$ be the minimum of $\delta_0$ in Lemma \ref{lemma_self} and all $\delta_{s_0,t_0}$ in Lemma \ref{lemma-tranverse} so that the conclusion for Lemma \ref{lemma_self} and \ref{lemma-tranverse} holds. Since there are only finitely many transverse intersection points {we have} $\delta>0$.  {The union over} all $t\in {\mathbb T}$ {of} $I_t = (t-\delta/4,t+\delta/4)$ covers ${\mathbb T}$. By compactness, we can find finitely many points $t_1,\cdots t_N$ so that ${\mathbb T}$ is covered by {the union of} $I_j = (t_j-\delta/4,t_j+\delta/4)$, {for} $j=1,\cdots, N$. 

     {We introduce a partition of unity subordinate to this covering.} Namely, we can find smooth functions $\varphi_j{\ge 0}$, for $j = 1,\cdots,N$, such that $\varphi_j$ is supported on $I_j$ and 
     $$
     \sum_{j=1}^N \varphi_j(t) = 1, \forall t\in {\mathbb T}.
     $$
Using (\ref{eq-auto-correlation}) and (\ref{eq-Gamma(s,t)}) and the partition of unity {we obtain}
\begin{equation}\label{eq-partition-unity}
\int f(x)d(\mu\ast\wt{\mu})(x) = \sum_{j=1}^N\sum_{k=1}^N\int_{I_j}\int_{I_k} f(\Gamma(s,t))~\varphi_j(s)\varphi_k(t)~dtds.
\end{equation}
Let $\varepsilon'>0$ be such that if $\gamma(\overline{I_j})\cap \gamma(\overline{I_k}) = \varnothing$, then the distance between the compact sets $\gamma(\overline{I_j}), \gamma(\overline{I_k})$ {satisfies}
$$
\mbox{dist} (\gamma(\overline{I_j}), \gamma(\overline{I_k})) >\varepsilon'>0.
$$
This $\varepsilon'>0$ guarantess that if ${\bf u}\in \Gamma(I_j\times I_k)$ and $\|{\bf u}\|<\varepsilon''$, then $\gamma(\overline{I_j})\cap \gamma(\overline{I_k}) \ne  \varnothing$. 

We now define $\varepsilon$ be the minimum of all $\varepsilon_{s_0,t_0}$ in Lemma {\ref{lemma-tranverse}},  {$\varepsilon_1$} defined in Lemma \ref{lemma-ball cover} and $\varepsilon'$ just defined.

Let  ${\bf u}\in B(0,\varepsilon)\setminus \{0\}$. Lemma \ref{lemma-ball cover} shows that it is possible to find $s,t$ such that $\Gamma (s,t) = {\bf u}$. We now claim that $\mu\ast\wt{\mu}$ has a smooth density function at $u$. This requires us to show that there exists a smooth function $g$ supported in an open neighborhood of $u$ such that 
$$
\int f(x) d(\mu\ast\wt\mu)(x) = \int f(x) g(x)dx
$$
for all continuous functions $f$ supported in that neighborhood.  From (\ref{eq-partition-unity}) and our definition of ${\varepsilon}$, 
\begin{equation}\label{eq-partition-unity1}
\int f(x)d(\mu\ast\wt{\mu})(x) = \sum_{\{(j,k): {\bf u}\in \Gamma(I_j\times I_k)\}}\int_{I_j}\int_{I_k} f(\Gamma(s,t))~\varphi_j(s)\varphi_k(t)~dtds.
\end{equation}
It suffices to show that each of the integrals above can be written as $
\int f({\bf u}) g_{{j,k}}({\bf u})d{\bf u}$ for some smooth functions $g_{{j,k}}$. There will be three cases:
\begin{enumerate}
    \item[Case 1:] {\bf  Same interval.} $I_j = I_k$.
    \item [Case 2:] {\bf  Adjacent interval.} $I_j \cap I_k\ne \varnothing$.
    \item[Case 3:]  {\bf Far away interval.} $I_j \cap I_k= \varnothing$.
\end{enumerate}

\noindent{\bf Case 1:} 
Consider the function {$F: I_j\times I_j\times B(0,\varepsilon)\to \RR^2$}, defined by 
\begin{equation}\label{eq_F}
F(s,t,{\bf x}) = \Gamma(s,t)- {\bf x}.
\end{equation}
By Lemma \ref{lemma_self}, $\Gamma$ is injective on $I_j\times I_j$ removing the diagonal {for any ${\bf x}$}. Let $(s_0,t_0)$ be such that $\Gamma(s_0,t_0) = {\bf u}$. Then $F(s_0,t_0,{\bf u}) =0$.  Note that the Jacobian of $F$ with respect to $(s,t)$ is equal to $J_{\Gamma}$, which is not zero at $(s_0,t_0).$ By the implicit function theorem, there exists a smooth function $g = (g_1,g_2)$ from a neigborhood of ${\bf u}$, denoted by $U$,  to a neighborhood of $(s_0,t_0)$ such that $s = g_1({\bf x})$, $t = g_2({\bf x})$, $(s_0,t_0) = g({\bf u})$ and $F(g({\bf x}),{\bf x}) = 0$. By the change of variable formula, 
\begin{equation}\label{eq-change-variable}
\int_{I_j}\int_{I_j} f(\Gamma(s,t))~\varphi_j(s)\varphi_j(t)~dtds= \int f({\bf x})\cdot \left(\varphi_j (g_1({\bf x}))\varphi_j (g_2({\bf x})) |J_{\Gamma}(g_1({\bf x}),g_2({\bf x}))|^{-1}\right) ~d{\bf x}
\end{equation}
whenever $f$ is supported on $U$. The function  in the bracket provides the smooth density in a neighborhood of ${\bf u}$.

\noindent{\bf Case 2:} Suppose that ${\bf u}\in \Gamma(I_j\times I_k)$ and $I_j\cap I_k\ne\varnothing$. Then the length of of the interval $I= I_j\cup I_k$ is at most $\delta<\delta_1$. Lemma \ref{lemma_self} shows that $\Gamma$ is injective on $I_j\times I_k$ removing the diagonal points. We can apply the same argument as in Case 1 on $F: I_j\times I_k\times B(0,\varepsilon)\setminus \{0\} \to \RR^2 $ with the same $F$ in (\ref{eq_F}) to obtain a smooth function. 

\noindent{\bf Case 3:} Suppose that ${\bf u}\in \Gamma(I_j\times I_k)$ and $I_j\cap I_k=\varnothing$. Let us write ${\bf u}  = \Gamma (s_1,t_1)$.   By our choice of $\varepsilon<\varepsilon''$, there must be some $s_0\in I_j$ and $t_0\in I_k$ such that $\gamma(s_0) = \gamma(t_0).$ By our choice of $\delta<\delta_{s_0,t_0}$ $I_j\subset I_{s_0,t_0}$ and $I_k\subset J_{s_0,t_0}$, where $I_{s_0,t_0}$ and $J_{s_0,t_0}$ are defined in Lemma \ref{lemma-tranverse}. Hence, $J_{\Gamma}(s_1,t_1)\ne 0$.  By the implicit function theorem, there exists $g = (g_1,g_2)$ from a neighborhood of ${\bf u}$, denoted by $U$, to a neighborhood of $(s_1,t_1)$ such that 
$$
\Gamma (g({\bf x})) = {\bf x}.
$$
Hence, we have a formula as in (\ref{eq-change-variable}), which provides us a smooth density. 

Combining three cases and (\ref{eq-partition-unity1}), we obtain a smooth density at ${\bf u}\in B(0,\varepsilon)\setminus\{0\}$. Since ${\bf u}$ is arbitrary, we establish a smooth density function for all ${\bf u}\in B(0,\varepsilon)\setminus \{0\}$.
\end{proof}

The proof of Theorem \ref{theorem-smooth} uses a partition of unity to produce a local density formula. We can also use Area formula (Theorem \ref{area formula}) to obtain a global density function. However, this function has multiple preimages and it is unclear {how} to show that the density function is smooth.  With Theorem \ref{theorem-smooth}, we are now ready to complete the proof for Theorem \ref{theorem2} using the same strategy we mentioned in the introduction.

\medskip

\begin{proof}[Proof of Theorem \ref{theorem2}.] Theorem \ref{theorem-smooth} showed that $\mu\ast\wt\mu$ has a smooth density in a punctured neighborhood of $0$. Suppose that the arc length measure $\mu$ is tight-frame spectral. Then Theorem \ref{prop-non-spectral} implies that $\delta_{\Lambda}$ has a positive spectral gap. However, $\mu$ is a singular measure in $\RR^2$, so any spectra must have a zero spectral gap by Theorem \ref{prop-spectral-gap}. This is a contradiction. Hence, the proof is complete. 
\end{proof}

\section{Discussions and Open questions.}\label{s:problems}

\subsection{Closed smooth curves.}This paper answered several natural questions about the spectrality of boundary of polygons and some other more general cases, but we did not completely classify all piecewise smooth curves that are spectral.  Concerning closed curves, it is natural to conjecture that 

\begin{conjecture}\label{conjecture2}
The arc length measure of any piecewise smooth closed curve is not spectral. 
\end{conjecture}
 Intuitively, closed curves would not be able to ``tile'' the space, so it should not be spectral (if one believes that Fuglede's conjecture is a guiding principle for determining the spectrality). In a simlar spirit, A domain with a ``hole" inside is known to be non-spectral \cite[Theorem 3.5]{lev2022fuglede}. We can attempt to prove this conjecture via a similar strategy in this paper. It is not hard to see that the support of $\mu\ast\wt\mu$ covers a neighborhood of the origin.  However, the existence of points with zero curvature or tangential self-intersections will destroy the smoothness of the density $\mu\ast\wt\mu$, while a flat line segment will create a singular part for the measure $\mu\ast\wt\mu$.  Theorem \ref{theorem1} and \ref{theorem2} deal with two extreme cases for Conjecture \ref{conjecture2}. We anticipate that proving the conjecture  will require a delicate study to interpolate the two scenarios. On the other hand, if Conjecture \ref{conjecture1} is correct, Conjecture \ref{conjecture2} is also automatically true. 

\subsection{Line spectra only?.} There are spectral measures supported on other non-closed curves. As mentioned in the introduction, the measure $\mu$ in (\ref{eq-symmetric-measure}) with $t = 0$, called the ``L-space"  in \cite{lai2021spectral}, is spectral. Indeed, when projecting to the subspace $y = -x$, the measure projects to a Lebesgue measure of a interval of length 2. From there, we obtain a spectrum $\{(n/2,-n/2): n\in\ZZ\}$.  This observation was generalized in  \cite[Theorem 1.2]{kolountzakis2025spectralitymeasureconsistingline}. In the same paper, it was also observed that all spectral measures supported on curves, to the best of our knowledge, are arising from a line spectrum in the same fashion as the $L$-space. This leads us naturally to the following question:

\begin{question}\label{line-spectra-only}
Is there a spectral curve in the plane (or collection of curves), under the arc length measure, that does not have a line spectrum?
\end{question}

Since an absolutely continuous spectral measure must have a constant density \cite{dutkay2014uniformity}, the question also asked for a spectral measure such that none of its projections {are a constant multiple of Lebesgue measure on a set that tiles the line.}

We further illustrate this question with two examples.

\begin{example}
    \rm{First, consider the arc-length measure on the two line segments connecting $(0, 0)$ to $(1, 0)$ and $(0, 1)$ to $(1, 1)$. This measure can be expressed as $\one_{[0,1]}dx\times \delta_{\{0,1\}}$. Hence, it has  a spectrum $\Lambda = \ZZ\times\Set{0, \frac12}$. See Fig.\ \ref{fig:stack-1}.}

\begin{figure}[h]
\ifdefined\SMART\resizebox{4cm}{!}{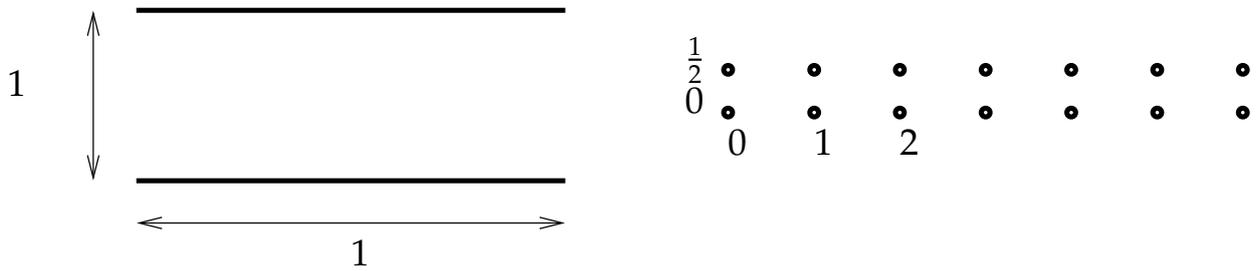}\else
\input stack-1.pdf_t
\fi
\caption{The arc-length measure $\mu$ on the two line segments shown left has the set $\ZZ\times\Set{0, \frac12}$ show on the right as a spectrum}\label{fig:stack-1}
\end{figure}
\textrm{It is however also true that $\mu$ has many line spectra. These are necessarily obtained by appropriate lines onto which $\mu$ projects to a spectral measure. There are many such lines (see Fig.\ \ref{fig:line-spectra}).}

\textrm{For instance, we can project onto the line perpendicular to the line joining points $(0, 1)$ and $(1, 0)$, in which case the projected measure is a constant on an interval. This is of course spectral and one spectrum of this interval (and, therefore of $\mu$ as well) is an arithmetic progression on the line of projection of spacing the reciprocal of the projected interval's length. }

\textrm{There are more spectra, which are line spectra that can be arbitrarily sparse in all directions (arbitrary sparseness of spectra is well-known for self-affine fractal measures {\cite{AL2023}}). Tilting the line of projection we can achieve that the projected measure will be a constant on the union of two intervals of the same length $\delta$ (this is always the case) and their gap is an integral multiple of $\delta$. Such a set of two intervals tiles the line and is therefore \cite{laba2001twointervals} also spectral in the line. Its spectrum, contained again in the line of projection, will again also be a spectrum of $\mu$. Notice though that the density of the spectrum is equal to $2\delta$ (the same as the projected set's measure) which can be arbitrarily small as $\delta$ can take arbitrarily small values.} 

\begin{figure}[h]
\ifdefined\SMART\resizebox{4cm}{!}{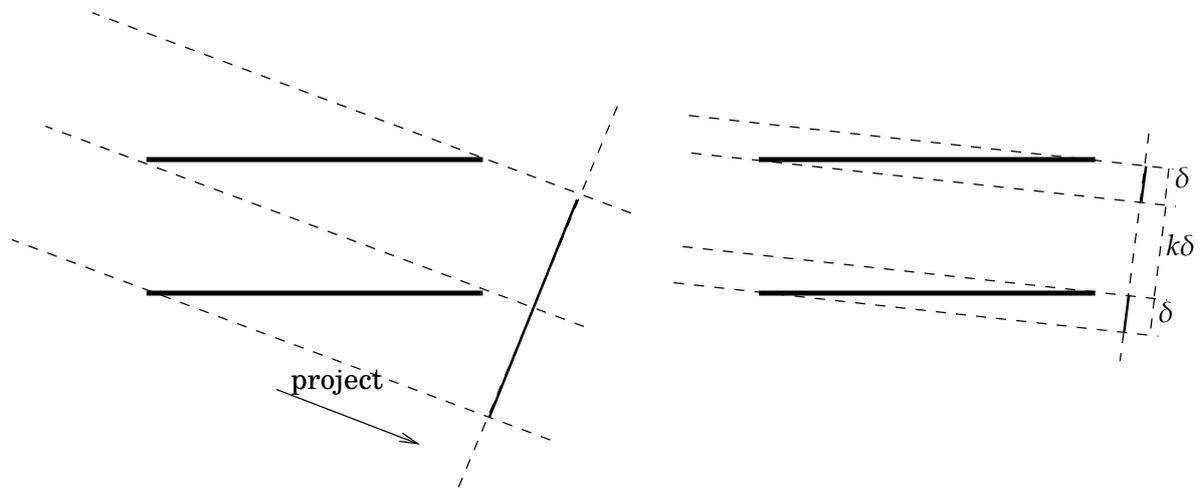}\else
\input line-spectra.pdf_t
\fi
\caption{The arc-length measure $\mu$ on the two line segments is projected, on the left, onto a single line segment, which gives a spectrum along the direction onto which we project. On the right, the same two line segments project to another spectral measure, as long as the gap of the projection is an integer multiple of the projected intervals.}\label{fig:line-spectra}
\end{figure}
\end{example}

\begin{figure}[h]
\ifdefined\SMART\resizebox{4cm}{!}{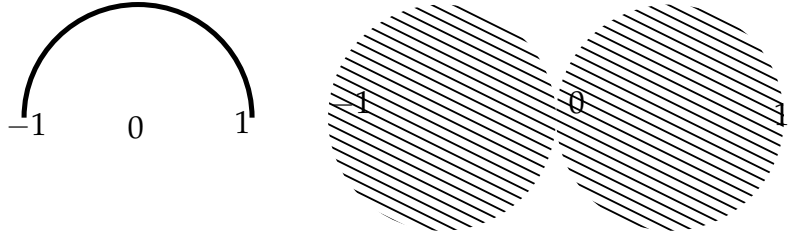}\else
\input semicircle.pdf_t
\fi
\caption{The arc-length measure $\mu$ on a semicircle gives rise to a measure $\mu*\wt\mu$ which is supported on the shaded region on the right. Thus it does not cover a neighborhood of the origin and our method is not applicable.}\label{fig:semicircle}
\end{figure}

\begin{example}
 {\rm Let us now consider the semi-circle of radius 1.  There is a spectral measure supported on the semi-circle. Namely, the push-forward of the Lebesgue measure of $[-1,1]$:
$$
\int f~d\nu = \frac12\int_{-1}^{1} f(x,\sqrt{1-x^2})~dx
$$
It is easy to see that $\frac12\ZZ\times \{0\}$ is a spectrum for $\nu$, which is a spectrum inside a line. However, the natural arc length measure $\mu$ is defined as follows:
$$
\int f~d\mu = \frac1{\pi}\int_{0}^{\pi} f(\cos x,\sin x)~dx.
$$
Note that since the projection of $\mu$ onto any subspaces are no longer Lebesgue measure with constant densities, there is no line spectrum. However, there may still be a chance for other spectra. Therefore, we ask:
\begin{question}\label{question2}
Is the arc length measure on the semi-circle a spectral measure?
\end{question}
We remark that the  support of $\mu\ast\wt\mu$ does not cover the origin in this case (see Fig.\ \ref{fig:semicircle}), so our method in this paper is not directly applicable. }
\end{example}

\subsection{Riesz bases of exponentials.} In \cite{iosevich2022fourier}, we know that all boundary of polygons (more generally polytopes in $\RR^d$) admits {a} frame of exponentials.  In this paper, we show however that  {boundaries} of polygons  are all non-spectral. Frames are overcomplete in general, an exact frame/ Riesz basis  is the frame that cease to be complete when one element is removed.  The natural question is this regard is now 
\begin{question}\label{question3}
Does the boundary of a polygon admit a Riesz bases of exponentials? 
\end{question}
 A partial result was obtained in \cite{lai2023riesz}, in which the boundary of square does not admit a type of structured exponentials as Riesz bases. 

\subsection{Fractal spectral measures.} Finally, our paper  sheds some light {into} fractal spectral measures. Conjecture \ref{conjecture1} in the introduction can provide some far reaching consequence in fractal spectral measures. As a simple example, we can prove the following.

\begin{proposition}
    Suppose that Conjecture \ref{conjecture1} holds. Let $\mu$ be a singular measure whose support has a positive Lebesgue measure. Then $\mu$ cannot be a spectral measure. 
\end{proposition}

\begin{proof}
    Suppose that $\mu$ is spectral with a spectrum $\Lambda$. Let $K$ be the support of $\mu$, which has positive Lebesgue measure. Notice that the support $\mu\ast\wt\mu$ is $K-K$, which contains an open set around the origin by the well-known Steinhaus theorem. If Conjecture \ref{conjecture1} holds, it means that ${\delta_{\Lambda}}$ has a positive spectral gap. However, this is a contradiction to Theorem \ref{prop-spectral-gap}.
\end{proof}

A typical example of singular measures without atoms whose support is $[0,1]$ is the Bernoulli convolution associated with the golden ratio.  It was not an easy proof in classifying the contraction ratios for which the Bernoulli convolution is a spectral measure \cite{HL2008,dai2012does}.  Conjecture \ref{conjecture1} leads to a much more natural and general approach to these problems. 

\noindent {\bf Acknowledgments.}  Chun-Kit Lai is partially supported by the AMS-Simons Research Enhancement Grants for Primarily Undergraduate Institution (PUI) Faculty.

\bibliographystyle{alpha}
\bibliography{mk-bibliography.bib}

\end{document}